\documentclass[a4paper,12pt]{article}

\usepackage[cp1251]{inputenc}
\usepackage[russian]{babel}

\usepackage{amsmath}
\usepackage{amssymb}
\usepackage{amsfonts}
\usepackage{amsbsy}
\newcommand{\lf}{\left}
\newcommand{\rt}{\right}
\usepackage{hyperref}

\begin{document}

УДК 519.233.5 + 519.654

\vspace{1.5cm}
\begin{center}
{\Large О связи между коэффициентами линейных регрессионных моделей разной размерности

\bigskip
В.Г. Панов}
\end{center}

\bigskip

Исследование регрессионных моделей различного типа является одним из основных методов исследования зависимости непрерывных переменных, несмотря на то, что сам метод регрессионного анализа является одним из самых старых в арсенале теоретической  и прикладной статистики \cite{Rao}--\cite{Fox}. Идеи, предложенные Ф.Гальтоном для решения интересовавшей его задачи \cite{DS, FR}, в дальнейшем получили глубокое и разнообразное развитие, приведшее к созданию регрессионного анализа как самостоятельной математической дисциплины. При этом потенциал развития этой теории далеко не исчерпан, и остаются широкие возможности для открытия новых подходов, развития, обобщения и углубления имеющихся методов и концепций как в самой теории регрессионного анализа, так и в его многочисленных приложениях \cite{Miller}--\cite{Harrell}.

Основным методом решения вариационной задачи, которая формулируется в регрессионном анализе, является, как правило, метод наименьших квадратов \cite{Linnik62}, возможно, с необходимыми обобщениями, продиктованными особенностями поставленной задачи. При этом сила и гибкость этого метода такова, что даже там, где он, кажется, не должен применяться, --- например, в логистической регрессии, ---  при надлежащей адаптации и обобщении он оказывается вполне состоятельным, наряду с методом наибольшего правдоподобия \cite{HL}--\cite{KK}.

При использовании регрессионного анализа в практической работе в распоряжении исследователя часто имеется только одна выборка, которую надо анализировать всеми возможными методами для получения максимальной доступной информации. Одна из типичных проблем прикладного регрессионного анализа состоит в выборе подходящего множества предикторов \cite{DS, Seber, Miller}. В данном случае для нас неважно, что будет критерием оценки качества модели с данным множеством предикторов. Мы будем рассматривать только линейные модели происхождение и значение которых для исследователя не входит в круг обсуждаемых вопросов. Существенно лишь то, что при изменении множества объясняющих переменных изменяется и сама регрессионная модель. Следовательно, возникает задача определения наличия и характера связи между моделями с разным набором предикторов или доказательство того, что этой связи нет.

Рассмотрим точную постановку задачи. Предварительно заметим, что ниже везде используется статистическая терминология (случайные переменные, математическое ожидание и т.д.), хотя по существу речь идет об аппроксимации данных линейными моделями с различным числом переменных (предикторов), и нахождение алгебраических связей между коэффициентами этих моделей. Таким образом, рассматриваемая задача имеет больше алгебраический характер, чем вероятностно-статистический. Однако, представляется, что наиболее существенной областью приложения доказанной теоремы будет именно регрессионный анализ, как в теоретической  так и в прикладной части. Да и сама задача более естественна именно в статистическом контексте.

Пусть заданы случайные переменные $X_1, X_2,\dots, X_n,Y,$ отличные от константы. В рассматриваемых ниже линейных регрессионных моделях случайные величины $X_1, X_2,\dots, X_n$ являются факторами, с помощью которых будет описываться изменение переменной-отклика $Y.$ В качестве уравнения, описывающего статистическую зависимость отклика $Y$ от предикторов $X_1, X_2,\dots, X_n$ рассмотрим уравнение множественной линейной регрессии
\begin{equation}
 \label{eq_Y_X_all}
y=b_0+\sum_{i=1}^{n}b_i\,x_i
\end{equation}

Постоянные $b_0,\ b_1,\dots,b_n$ являются решениями вариационной задачи минимизации среднеквадратичного уклонения (\textsf{E} ---
оператор математического ожидания)

\begin{equation}\label{task_b}
\min\limits_{b_0,b_1,\dots,b_n}\text{\textsf{E}}\lf(Y-b_0-\sum_{k=1}^nb_kX_k\rt)^2
\end{equation}

Во множестве объясняющих переменных переменных $X_1,\ X_2,\dots,\ X_n$ зафиксируем некоторое подмножество
$\{X_{i_1},X_{i_2},\dots,X_{i_k}\}.$ Так как, вообще говоря, между переменными  $X_1,\ X_2,\dots,\ X_n$ возможны  различные статистические связи, то имеет смысл рассмотреть уравнения множественной регрессии между самими предикторами, а именно, линейные уравнения, выражающие зависимость каждой переменной  $X_i$ от выбранной системы объясняющих переменных $\{X_{i_1},X_{i_2},\dots,X_{i_k}\}$ (все коэффициенты этих уравнений зависят также от множества переменных $\{X_{i_1},X_{i_2},\dots,X_{i_k}\}$)
\begin{equation}
 \label{eq_X_X_ij}
x_i=c_{i}+\sum_{j=1}^{k}c_{i,i_j}\,x_{i_j},\qquad i=1,2,\dots,n,
\end{equation}
где коэффициенты $c_{0ij},\ c_{ij}$ являются решениями вариационной задачи

\begin{equation}\label{task_c}
\min\limits_{c_i,c_{i,i_j}}\text{\textsf{E}}\lf(X_i-c_{i}-\sum_{j=1}^{k}c_{i,i_j}\,X_{i_j}\rt)^2,\quad  i,j=1,2,\dots,n
\end{equation}

При этом будем считать, что если $i\in \{i_1,i_2,\dots,i_k\},$ то $c_{i}=0,\ c_{i,i_j}=\delta_{ii_j}.$

Наконец, множественной линейной регрессией будем также описывать зависимость переменной $Y$ от предикторов
$\{X_{i_1},X_{i_2},\dots,X_{i_k}\}$
\begin{equation}
 \label{eq_Y_X_ij}
y=a_{0}+\sum_{j=1}^{k}a_{i_j}\,x_{i_j},
\end{equation}
где коэффициенты $a_0,\ a_{i_j}$, как и для других вариационных задач, являются решениями задачи минимизации среднеквадратичного
уклонения

\begin{equation}
\min\text{\textsf{E}}\lf(Y-a_{0}-\sum_{j=1}^ka_{i_j}X_{i_j}\rt)^2 \label{task_a}
\end{equation}

При обычных предположениях (см., например, \cite{Linnik62, MagnusNeu, KozlovProchorov}) каждая из этих задач имеет единственное решение, которое находится из соответствующей системы линейных уравнений.

Рассматриваемая ниже задача состоит в исследовании связи между коэффициентами введенных выше регрессионных моделей, т.е. существует ли какая-либо связь между коэффициентами $\{a_{i_j}\},\{b_k\},\{c_{i,i_j}\}$? В частном случае простых (парных) регрессионных уравнений между предикторами и между откликом и каждым из предикторов эта задача была исследована в \cite{PV}.

Проблема изучения связи регрессионных моделей разной размерности возникла из задачи анализа статистической зависимости признаков при согласованном изменении предикторов в работах \cite{VMChA, MVCh}. Там же была обоснована необходимость изучения такого рода задач в связи с проблемами медико-биологического мониторинга. В работе \cite{MVCh} были предложены два варианта согласованного изменения предикторов $\{X_1,X_2,\dots,X_n\}$: пропорциональное и корреляционное. В данной работе изучаются соотношения между коэффициентами множественных регрессионных моделей \eqref{eq_Y_X_all}, \eqref{eq_X_X_ij},\eqref{eq_Y_X_ij} при множественной регрессионной зависимости предикторов \eqref{eq_X_X_ij}.

Отметим, некоторые особенности поставленной задачи. Во-первых, мы не рассматриваем задачу статистического оценивания параметров линейной регрессионной модели и связи между статистическими оценками параметров разных моделей. В данном случае предметом исследования является связь между регрессионными коэффициентами моделей с разным множеством предикторов, которые были вычислены по одному и тому же множеству данных. Следовательно, здесь не возникает вариативность выборочных данных, но по одной и той же выборке появляется вариативность линейных моделей, зависящих от множества выбираем для описания предикторов. Во-вторых, не затрагивается хорошо изученный вопрос качества (проверки гипотез) линейных моделей \cite{DS, Seber, Miller, LL}. Регрессионные модели, которые были введены выше, могли появиться из самых разных соображений, не обязательно продиктованных каким-либо критерием качества подгонки. В-третьих, сформулированная ниже теорема доказана только для того случая, когда коэффициенты регрессионных моделей найдены обычным (классическим) методом наименьших квадратов. Вследствие того, что все вычисления делаются по одному  и тому же массиву данных,  применение МНК в данном случае носит не статистический, а аппроксимационный характер. Таким образом, можно сказать, что доказанная теорема выражает свойство самого метода наименьших квадратов.

Задача определения связи между коэффициентами линейных регрессионных моделей разной размерности, насколько смог установить автор,  в изложенной выше постановке ранее не рассматривалась. По-видимому, причина этого в том, что никакой связи между коэффициентами различных моделей не предполагалось. Наиболее близко к рассматриваемой задаче стоит задача расширения множества предикторов, рассмотренная впервые в \cite{Cochran} и обобщенная в \cite{Quenouille}. Изложение этих методов дано в \cite{Seber}. Здесь важно отметить, что в работах \cite{Cochran, Quenouille} рассматривается задача определения коэффициентов уравнения линейной регрессии при добавлении некоторого множества предикторов только для добавляемых предикторов, т.е. предикторы при старых переменных сохраняются теми же, какими они были без добавленных переменных. Таким образом, новая модель <<помнит>> свое происхождение и поэтому не является чисто новой. В рассматриваемой ниже задаче никакой исходной связи между моделями с различным множеством предикторов не предполагается, и коэффициенты этих моделей вычисляются независимо друг от друга с помощью метода наименьших квадратов.

Кроме того, в регрессионном анализе хорошо известен метод  построения уравнения множественной регрессии путем последовательного включения в уравнение каждой новой переменной \cite[раздел 22.7]{DS}. Существенную роль в этом подходе играет вычисление регрессии остатков на остатки, которое приводит к значительному увеличению объема вычислений и накоплению вычислительной ошибки. В данном случае этого не будет требоваться, так как коэффициенты уравнения множественной регрессии отклика можно получить путем матричного умножения через коэффициенты парных регрессионных уравнений, если в качестве исходных регрессионных моделей \eqref{eq_X_X_ij} выбрать всевозможные уравнения простых регрессий предикторов, а для модели  \eqref{eq_Y_X_ij} (для каждого предиктора) взять также уравнение простой регрессии отклика по данному предиктору. В работе \cite{PV} это обсуждается подробнее.

Перед формулировкой основного утверждения докажем  вспомогательные леммы. Для случайной величины $X$ обозначим $\overline{X}$ среднее
значение (математическое ожидание) этой случайной величины. Так как мы будем рассматривать только дискретные случайные величины, то ниже во всех случаях среднее значение соответствующей случайной величины существует.

В соответствии с методом наименьших квадратов \cite{Linnik62, MagnusNeu} решение оптимизационных задач \eqref{task_b}, \eqref{task_c}, \eqref{task_a} сводится к решению следующих систем линейных уравнений ($\mathbf{X},\mathbf{X}_I$ --- матрицы систем линейных уравнений, соответствующих задачам  \eqref{task_b}, \eqref{task_c}, \eqref{task_a}, \ $\mathbf{YX},\mathbf{XX}_I,\mathbf{YX}_I$ --- столбцы свободных членов этих систем)

\begin{equation}\label{mnq_b}
\mathbf{B}\cdot\mathbf{X}=\mathbf{YX}
\end{equation}

\begin{equation}
\label{mnq_c}
\mathbf{C}_I\cdot\mathbf{X}_I=\mathbf{XX}_I
\end{equation}

\begin{equation}
\label{mnq_a}
\mathbf{A}_I\cdot\mathbf{X}_I=\mathbf{YX}_I,
\end{equation}
где через $I$ обозначено множество индексов выделенной системы предикторов \( I=\{i_1,i_2,\dots,i_k\}\) и введены следующие матрицы
\[
\mathbf{A}_I=(a_{i_1},a_{i_2},\dots,a_{i_k})_{1\times k};
\ \mathbf{B}=(b_1,b_2,\dots,b_n)_{1\times n};
\ \mathbf{C}_I=(c_{i,i_j})_{n\times k},
\]
причем $c_{i,i_j}=\delta_{ii_j},$ если $i\in I.$

Зафиксируем индекс $i_j\in I$ и рассмотрим следующую квадратную матрицу порядка $n+k+2$ ($\mathbf{0}$ означает нулевую матрицу
размерности $(n+1)\times k$)

$$
\mathcal{A}=\left(
  \begin{array}{cc}
    M & \mathbf{0} \\
    N & N_{i_j} \\
  \end{array}
\right),
$$
где $M$ --- расширенная матрица системы \eqref{mnq_b} (матрица размерности $(n+1)\times(n+2)$);
$N$ --- матрица размерности $(k+1)\times(n+2)$
$$
N=\left(
  \begin{array}{cccccc}
    0 & \overline X_1 & \overline X_2 & \cdots& \overline X_n & \overline Y \\
    0 & \overline {X_{i_1}X_1} &  \overline {X_{i_1}X}{\hspace{-1pt}}_2 & \cdots &  \overline {X_{i_1}X}{\hspace{-1pt}}_n &  \overline
    {X_{i_1}Y} \\
    0 &  \overline {X_{i_2}X}{\hspace{-1pt}}_1 &  \overline {X_{i_2}X}{\hspace{-1pt}}_2 &\cdots &  \overline
    {X_{i_2}X}{\hspace{-1pt}}_n &  \overline {X_{i_2}Y} \\
    \vdots &  \vdots & \vdots & \cdots & \vdots & \vdots \\
    0 &  \overline {X_{i_k}X}{\hspace{-1pt}}_1 &  \overline {X_{i_k}X}{\hspace{-1pt}}_2 & \cdots &  \overline
    {X_{i_k}X}{\hspace{-1pt}}_n &  \overline {X_{i_k}Y} \\
  \end{array}
\right);
$$
матрица $N_{i_j}$ имеет размерность $(k+1)\times k$ и получается из $\mathbf{X}_I$ (матрицы систем \eqref{mnq_c} и \eqref{mnq_a})
вычеркиванием столбца
$$
\left(
  \begin{array}{c}
    \overline{X}{\hspace{-1pt}}_{i_j}\\
    \overline{X_{i_1}X}{\hspace{-1pt}}_{i_j} \\
    \dots \\
    \overline{X_{i_k}X}{\hspace{-1pt}}_{i_j} \\
  \end{array}
\right)
$$

\textbf{Лемма 1. }Матрица $\mathcal{A}$ вырождена.

\textbf{Доказательство.} Изменением нумерации предикторных переменных доказательство леммы можно свести к случаю $I=\{1,2,\dots,k\},\
i_j=k.$ Вычитая из $i$-ой строки этого определителя его $(n+i+1)$-ю строку, $i=1,2,\dots,k+1$, а затем вычитая из первого столбца $(n+3)$-й
столбец, получим (индексами указаны
размерности соответствующих матриц)

\begin{equation}\label{detA_n}
\det(\mathcal{A})=\left|
                \begin{array}{cc}
                  \mathbf 0_{(k+1)\times(n+2)} & -N_{(k+1)\times k} \\
                  P_{(n-k)\times(n+2)} & \mathbf 0_{(n-k)\times k} \\
                  Q_{(k+1)\times(n+2)} & N_{(k+1)\times k} \\
                \end{array}
              \right|
\end{equation}

Теперь, считая $k$ фиксированным, проведем индукцию по $n$. Минимально допустимое значение для $n$ равно $k+1,$ или $k=n-1.$ Тогда

\begin{equation}\label{det_A_n=k+1}
\det(\mathcal A)= \left|
                    \begin{array}{cc}
                    \mathbf 0_{n\times(n+2)} & -N_{n\times (n-1)} \\
                  P_{1\times(n+2)} & \mathbf 0_{1\times (n-1)} \\
                  Q_{n\times(n+2)} & N_{n\times (n-1)} \\
                    \end{array}
                  \right|
\end{equation}

Раскрывая этот определитель по строке $P_{1\times(n+2)},$ получим линейную комбинацию определителей порядка $2n$, имеющих вид
$$
\left|
  \begin{array}{cc}
    \mathbf 0_{n\times(n+1)} & -N_{n\times (n-1)} \\
    Q'_{n\times(n+1)} & N_{n\times (n-1)} \\
  \end{array}
\right|
$$

Сгруппируем элементы этого определителя, присоединив $(n+1)$-й столбец к остальным $n-1$ столбцам. В результате получим представление
определителя в виде блочной матрицы

\[  \left|
  \begin{array}{c|c}
    \mathbf 0_{n\times n}  &\mathbf{0}_{n\times1}\  {-N_{n\times (n-1)}} \\
    \hline
    Q''_{n\times n} & N'_{n\times n} \\
  \end{array}
\right|
\]

В двух последних определителях через $Q'$ и $Q''$ обозначены матрицы, которые получаются из матрицы $Q$ вычеркиванием одного или двух
столбцов соответственно. Аналогично
матрица $N'$ получается из матрицы $N$ добавлением одного столбца (вычеркнутого из матрицы $N$).

Так как определитель, стоящий в правом верхнем углу равен нулю, то и весь определитель \eqref{det_A_n=k+1} равен нулю, что доказывает
случай леммы для $n=k-1.$

Пусть лемма доказана для $n\geqslant k+2.$ Раскроем определитель \eqref{detA_n} по первой строке матрицы $P_{(n-k)\times(n+2)}$. Тогда
получим линейную
комбинацию определителей вида

\begin{equation}\label{detA_n-1}
\left|
                \begin{array}{cc}
                  \mathbf 0_{(k+1)\times(n+1)} & -N_{(k+1)\times k} \\
                  P_{(n-k-1)\times(n+1)} & \mathbf 0_{(n-k-1)\times k} \\
                  Q_{(k+1)\times(n+1)} & N_{(k+1)\times k} \\
                \end{array}
              \right|
\end{equation}

\medskip
Сопоставляя определители \eqref{detA_n} и \eqref{detA_n-1}, замечаем, что определитель \eqref{detA_n-1} имеет вид определителя
\eqref{detA_n} при замене $n$ на
$n-1$. Следовательно, по предположению индукции, все эти определители равны нулю. Лемма 1 доказана.

\medskip

Для формулировки и доказательства леммы 2 нам понадобится теорема Лапласа, которую удобно сформулировать в виде \cite{Marcus-Minc}
(см. также \cite{Gantmakher}).

Пусть $Q_{r,n}$ --- совокупность всех строго возрастающих последовательностей длины $r$ из $(1,2,\dots,n);$
$A\binom{\alpha}{\beta}=A\binom{\alpha_1,\alpha_2,\dots,\alpha_r}{\beta_1,\beta_2,\dots,\beta_r}$ --- определитель подматрицы матрицы
$A,$ образованной пересечением
строк с номерами из $\alpha=(\alpha_1,\alpha_2,\dots.\alpha_r)$ и столбцов с номерами из $\beta=(\beta_1\beta_2,\dots,\beta_r)$ где
$\alpha,\ \beta$ --- наборы из
$Q_{r,n}$; $A\binom{\alpha^{\ast}}{\beta^{\ast}}$ --- определитель подматрицы матрицы $A,$ стоящей на пересечении строк с номерами из
$\alpha^{\ast}$ и столбцов с
номерами из $\beta^{\ast},$ где $\alpha^{\ast},\beta^{\ast}$ --- наборы из $Q_{n-r,n},$ дополнительные для наборов $\alpha$ и $\beta$
соответственно:
$\alpha\cap\alpha^{\ast}=\varnothing,\ \beta\cap\beta^{\ast}=\varnothing,\
\alpha\cup\alpha^{\ast}=\beta\cup\beta^{\ast}=\{1,2,\dots,n\};$ для любого набора $\beta =
(\beta_1,\beta_2,\dots,\beta_r)$ обозначим $|\beta|=\sum\limits_{i=1}^r\beta_i.$

\textbf{Теорема Лапласа}. Пусть $A$ --- квадратная матрица порядка $n,\ 1\leqslant r\leqslant n,\ \alpha$ --- фиксированная
последовательность из $Q_{r,n}$. Тогда
$$
\det(A)=(-1)^{|\alpha|}\sum_{\beta\in Q_{r,n}}(-1)^{|\beta|}A\binom{\alpha}{\beta}\cdot A\binom{\alpha^{\ast}}{\beta^{\ast}}
$$

Возьмем $r=n+1,\alpha =(1,2,\dots,n,n+1),\alpha\in Q_{n+1,n+k+2}$ и применим к матрице $\cal A$ теорему Лапласа. Учитывая лемму 1,
получим равенство

\begin{equation}
\label{Laplace_A} 0=\sum_{\beta\in Q_{n+1,n+k+2}} (-1)^{|\beta|}\cdot \mathcal{A}\binom{1,2,\dots,n,n+1}\beta \cdot
\mathcal{A}\binom{n+2,\dots,n+k+2}{\beta^{\ast}},
\end{equation}
где $\beta^{\ast}$ --- набор, дополнительный для набора $\beta$ в множестве $(1,2,\dots,n+k+2),$ т.\,е.
$\beta\cup\beta^{\ast}=(1,2,\dots,n+k+2).$

\textbf{Лемма 2.} Количество наборов $\beta,$ для которых произведение определителей в сумме \eqref{Laplace_A} отлично от нуля, не
более $n+1.$

\textbf{Доказательство.} Пусть $\beta$ --- некоторый набор из $Q_{n+1,n+k+2}.$ Если 1 не входит в $\beta,$ то дополнительный набор
$\beta^{\ast}$ начинается с
1 и, следовательно,
\[
\mathcal{A}\binom{n+2,\dots,n+k+2}{1\qquad\qquad\ast\quad}=0,
\]
так как в силу строения матрицы $\cal A$ определитель $\mathcal{A}\binom{n+2,\dots,n+k+2}{1\quad\quad\ast\quad}$ содержит нулевой
столбец (звездочкой
обозначены все числа дополнительного набора $\beta^{\ast},$ кроме 1).

С другой стороны, если в набор $\beta$ входит любое из чисел $n+3,\ n+4,\dots,n+k+2,$ то
\[
\mathcal{A}\binom{n+2,\dots,n+k+2}{\ast\quad n+j+2\quad\ast\quad}=0,\qquad j=1,2,\dots,k
\]

Следовательно, в сумме \eqref{Laplace_A} ненулевыми могут быть только слагаемые, соответствующие наборам $\beta\in Q_{n+1,n+k+2}$,
которые начинаются с 1 и
заканчиваются не более чем на $n+2.$

Таким образом, интересующие нас наборы $\beta$ должны иметь вид
\[
(1,\beta_2,\ldots,\beta_n,n+1) \qquad \text{или} \qquad (1,\beta_2,\ldots,\beta_n,n+2)
\]

Так как первый набор определен однозначно и равен $(1,2,\dots,n,n+1)$, то остается доказать, что наборов второго типа будет ровно $n.$
Из этих наборов выделим
случай
$$
\beta = (1,2,\ldots,n,n+2)
$$
Для остальных наборов рассматриваемого типа существует единственное число $i$ от $2$ до $n,$ (включительно) такое, что
\begin{equation} \label{beta_i}
\begin{aligned}
\beta_j&=j\quad \text{ при } j=1,2,\dots,i\\
\beta_i&=i+1,\text{ при } j=i  \\
\beta_j&=j+1\quad\text{при } j=i+1, i+2,\ldots,n+1
\end{aligned}
\end{equation}
Обозначим через $\mathfrak{b}_i$ последовательность $\beta\in Q_{n+1,n+k+2}$, определяемую условиями \eqref{beta_i}. Таким образом,
наборов \eqref{beta_i} будет ровно
столько, сколько значений может принимать индекс $i$ в этих равенствах. Следовательно, всего таких наборов будет $n-1+1=n,$ что и
доказывает лемму.

Как нетрудно подсчитать, для наборов $\beta$ из леммы 2 имеют место равенства:
\[
\begin{aligned}
|(1,2,\ldots,n,n+1)|&=\frac{(n+1)(n+2)}2\\
|(1,2,\ldots,n,n+2)|&=\frac{(n+1)(n+2)}2+1\\
|(1,2,\ldots,i-1,i+1,\ldots,n+1,n+2)|&=\frac{(n+2)(n+3)}2-i
\end{aligned}
\]

Теперь сформулируем и докажем основное утверждение.

\bigskip
\textbf{Теорема. }Если $\{a_{i_j}\},\ \{b_{i}\},\ \{c_{i,i_j}\}$ --- решения задач \eqref{task_b}, \eqref{task_c}, \eqref{task_a}, то
имеет место равенство

\begin{equation}\label{main_eq}
a_{i_j}=\sum_{i=1}^{n}b_ic_{i,i_j}
\end{equation}

\textbf{Доказательство. }Введем следующий определитель $(k+1)$-го порядка (главный определитель систем \eqref{mnq_c} и \eqref{mnq_a})
$$
\delta=\mathrm{Det}(\mathbf{X}_I)
$$

Обозначим через $\delta_{i_j}$ определитель, который  получается из определителя $\delta$ заменой $(j+1)$-го столбца столбцом
$\mathbf{YX}_I$ из системы \eqref{mnq_a}, $j=1,\ 2,\dots,k$. При замене первого столбца определителя $\delta$ столбцом $\mathbf{YX}_I$ получим определитель, который обозначим $\delta_{0}.$ Определитель $d_{i,i_j}$ получается из определителя $\delta$ заменой $(j+1)$-го столбца столбцом $\mathbf{XX}_I$ из системы \eqref{mnq_c}, $i=1,2,\dots,n,\ j=1,\
2,\dots,k$. При замене первого столбца определителя $\delta$ столбцом свободных членов системы \eqref{mnq_c} получим определитель
$d_{i}.$ Главный определитель системы \eqref{mnq_b} обозначим $\Delta$, а определители $\Delta_i$ получаются из определителя $\Delta$ заменой $(i+1)$-го столбца столбцом $\mathbf{YX}$ из  системы \eqref{mnq_b}, $i=0,1,\dots,n.$

Мы будем считать, что определители $\delta$ и $\Delta$ не равны нулю, т.\,е. решения этих систем существуют и единственны. Тогда по
формулам Крамера имеют место
равенства
\begin{gather}
a_{0}=\dfrac{\delta_{0}}{\delta},\qquad a_{i_j}=\dfrac{\delta_{i_j}}{\delta}\label{solve_a}\\
c_{i}=\dfrac{d_{i}}{\delta},\qquad c_{i,i_j}=\dfrac{d_{i,i_j}}{\delta}\label{solve_c}\\
\label{solve_b} b_0=\dfrac{\Delta_0}{\Delta},\qquad b_i=\dfrac{\Delta_i}{\Delta}
\end{gather}

Рассмотрим определители из равенства \eqref{Laplace_A} с учетом леммы 2. Напомним, что из матрицы определителя ${\cal A}$, имеющей размерность $(k+1)\times(k+1)$ и расположенной в правом нижнем углу этого определителя, выброшен столбец
$$
\left(
  \begin{array}{c}
    \overline{X}{\hspace{-1pt}}_{i_j}\\
    \overline{X_{i_1}X}{\hspace{-1pt}}_{i_j} \\
    \dots \\
    \overline{X_{i_k}X}{\hspace{-1pt}}_{i_j} \\
  \end{array}
\right)
$$

Имеем следующие равенства

\begin{align*}
{\cal A}\binom{1,2,\ldots,n+1}{1,2,\ldots,n+1}&=\Delta, &{\cal A}\binom{n+2,\dots,n+k+2}{n+2,\dots,n+k+2} &=(-1)^j\delta_{i_j},\\
{\cal A}\binom{1,2,\ldots,n+1}{1,2,\ldots,n+2}&=\Delta_n, &{\cal A}\binom{n+2,\dots,n+k+2}{n+1,n+3,\dots,n+k+2}
&=(-1)^jd_{n,\,i_j},\\
{\cal A}\binom{1,2,\ldots,n+1}{\mathfrak{b}_i}&=(-1)^{n-i+1}\Delta_{i-1}, &{\cal A}\binom{n+2,\dots,n+k+2}{i,n+3,\dots,n+k+2}
&=(-1)^jd_{i-1\,,i_j},\\
\text{ где $\mathfrak{b}_i$ --- набор } \eqref{beta_i},&\  i=2,3,\ldots,n.\\
\end{align*}

Подставляя эти выражения в равенство \eqref{Laplace_A}, получаем
$$
\begin{aligned}
0=&\sum_{\beta\in Q_{n+1,n+k+2}} (-1)^{|\beta|}\cdot \mathcal{A}\binom{1,2,\dots,n,n+1}\beta \cdot
\mathcal{A}\binom{n+2,\dots,n+k+2}{\beta^{\ast}}=\\=
&(-1)^{\frac{(n+1)(n+2)}2}\Delta\cdot (-1)^j\delta_{i_j}+ (-1)^{\frac{(n+1)(n+2)}2+1}\Delta_n\cdot
(-1)^jd_{n,\,i_j}+\\+&\sum_{i=2}^n(-1)^{\frac{(n+2)(n+3)}2-i}(-1)^{n-i+1}\Delta_{i-1}\cdot (-1)^jd_{i-1,\,i_j}=\\=&
(-1)^{\frac{(n+1)(n+2)}2+j}\lf[\Delta\cdot
\delta_{i_j}-\Delta_n\cdot d_{n,\,i_j}-\sum_{i=2}^n(-1)^{n-i}(-1)^{n-i}\Delta_{i-1}\cdot
d_{i-1\,i_j}\rt]=\\=&(-1)^{\frac{(n+1)(n+2)}2+j}\lf[\Delta\cdot
\delta_{i_j}-\Delta_n\cdot d_{n,\,i_j}-\sum_{i=2}^n\Delta_{i-1}\cdot d_{i-1,\,i_j}\rt]
\end{aligned}
$$
Следовательно,
$$
\Delta\cdot \delta_{i_j}=\Delta_n\cdot d_{n,\,i_j}+\sum_{i=2}^n\Delta_{i-1}\cdot d_{i-1,\,i_j}
$$
Деля это равенство на $\delta\cdot \Delta,$ и учитывая равенства \eqref{solve_a}--\eqref{solve_b}, получим
$$
a_{i_j}=\sum_{i=1}^n b_i\cdot c_{i,\,i_j},
$$
что и требовалось доказать.

С использованием введенных перед леммой 1 обозначений эта теорема принимает вид
\[
\mathbf{A}_I=\mathbf{B}\cdot \mathbf{C}_I
\]
Это означает, что коэффициенты равенства \eqref{eq_Y_X_ij} получаются из равенства \eqref{eq_Y_X_all} формальной подстановкой правых
частей уравнений \eqref{eq_X_X_ij} с
последующим приведением подобных. Однако поскольку все коэффициенты в уравнениях \eqref{eq_Y_X_all}, \eqref{eq_X_X_ij},
\eqref{eq_Y_X_ij} находятся из условия
минимизации соответствующего квадратичного функционала \eqref{task_b}, \eqref{task_c}, \eqref{task_a}, то доказываемое равенство
\eqref{main_eq} не является следствием
только линейных уравнений \eqref{eq_Y_X_all}, \eqref{eq_X_X_ij}, \eqref{eq_Y_X_ij}. Можно сказать, что равенство \eqref{main_eq}
является свойством метода наименьших
квадратов, с помощью которого находятся коэффициенты  $a_{i_j}, b_k, c_{i,i_j}$. Это свойство можно сформулировать как условие
коммутирования нахождения решения задачи
\eqref{task_a} с помощью метода наименьших квадратов и матричного умножения.

Рассмотрим некоторые следствия доказанной теоремы.

\textbf{Следствие 1.} Если $I=\{i\},$ то
\[
a_i= b_i+\sum_{j=1,\ j\neq i}^nb_jc_{j\hspace{1pt}i}=\sum_{j=1}^nb_jc_{j\hspace{1pt}i},\qquad\text{где }c_{ii}=1.
\]

\textbf{Следствие 2. }\cite{PV}. Если коэффициенты $a_i,\ c_{ij}$ определены для любых $i, j$ (т.\,е. рассматриваются простые корреляции между
парами предикторов
$X_i,X_j$ и, соответственно, простые корреляции $Y$ по $X_i,i=1,2\dots,n$), то имеет место равенство
\[
A=B\cdot C,
\]
где $A=(a_i)_{1\times n},B=(b_j)_{1\times n},C=(c_{j\hspace{1pt}i})_{n\times n}.$

\textbf{Следствие 3. }\cite{PV}. Если в условиях следствия 2 матрица $C$ обратима, то
\[
B=A\cdot C^{-1}
\]

\textbf{Следствие 4. }Пусть $\{1,2,\dots,n\}=I_1\cup I_2\cup\dots\cup I_m$ --- объединение непересекающихся множеств. Пусть $A_{I_j}$
--- строка коэффициентов
$(a_{l}),l\in I_j.$ Обозначим через $\mathbf{A}$ строку коэффициентов, полученную последовательным приписыванием строк $A_{i_j}$
справа к строке $A_{i_1}$. Аналогично
получим матрицу $\mathbf{C}$ из матриц $C_{I_j}$ (последовательным присоединением матриц $C_{I_j}$ справа к матрице $C_{i_1}$). Тогда
\[
\mathbf{A}=B\cdot \mathbf{C}
\]

Коэффициенты $a_{j}$ можно определить не только для случая $j\in I,$ но и для всех $j=1,2,\dots, n.$ А именно, если $j\not\in I,$ то
по определению будем считать, что
$a_{j}=0.$ Это согласуется с тем, что в уравнении \eqref{eq_Y_X_ij} присутствуют только слагаемые с коэффициентами $a_{j},j\in I.$
Тогда строка $A$ всегда будет иметь
размерность $1\times n.$

Также можно естественно доопределить значения коэффициентов $c_{i,j}$ для любых $i,j=1,2,\dots,n $ следующим образом
\[
c_{i,j}=\begin{cases} 0,& \text{ если $j\not\in I;$}\\
\delta_{j},&\text{ если $i,j\in I;$}\\
c_{i,j},&\text{ если $i\not\in I,j\in I$}
\end{cases}
\]
В последнем случае равенство означает, что значение $c_{i,j}$ находится из системы \eqref{mnq_c}, т.\,е. является решением
соответствующей оптимизационной задачи
\eqref{task_c}.

Таким образом, матрица $C$ при любом множестве $I$ будет определена как матрица размерности $n\times n.$

\textbf{Следствие 5. } Если матрицы $A,C$ определены так, как это описано выше,  то для любого множества $I$ имеет место равенство
\[
A=B\cdot C
\]
где $B=(b_i)_{1\times n}$ --- матрица коэффициентов уравнения множественной регрессии \eqref{eq_Y_X_all}.

\bigskip
В заключение автор хотел бы выразить благодарность профессору А.Н. Вараксину за привлечение внимания к той замеченной им числовой
закономерности, из которой выросла настоящая работа, а также доценту Ю.В. Нагребецкой за высказанные замечания и плодотворные дискуссии по некоторым алгебраическим вопросам.


\begin{thebibliography}{99}

\bibitem{Rao} Рао С.Р. Линейные статистические методы и их приенение. -- М.: Наука, 1968.--548 с.
\bibitem{DS} Дрейпер Н., Смит Г. Прикладной регрессионный анализ. 3-е изд. -- Диалектика, 2007.--912 с.
\bibitem{Seber} Себер Дж. Линейный регрессионный анализ. -- М.: Мир, 1980.--456 с.
\bibitem{Gross} Gross J. Linear regression. -- Springer, 2003.--394 p.
\bibitem{Fox} Fox L. Applied Regression Analysis, Linear Models, and Related Methods. --- Sage, 1997.--597 p.
\bibitem{FR} Фёрстер Э., Рёнц Б. Методы корреляционного и регрессионного анализа. ---  М.: Финансы и статистика. 1983.-- 303 с.
\bibitem{Miller} Miller A. Subset selection in regression. 2nd ed. --- Chapman \& Hall/CRC, 2002.--238 p.
\bibitem{Berk} Berk R. A. Regression analysis: a constructive critique. - SAGE, 2004.--259 p.
\bibitem{BKW} Belsley D.A., Kuh E. Welsch R.E. Regression Diagnostics: Identifying Influential Data and Sources of Collinearity. --- Wiley-IEEE, 2004.--292 p.
\bibitem{Harrell} Harrell F.E. Regression modeling strategies: with applications to linear models, logistic regression, and survival analysis. --- Springer, 2001.--568 p.
\bibitem{Linnik62} Линник Ю.\,В. Метод наименьших квадратов и основы математико-статистической теории
обработки наблюдений. 2-е изд. М.: Наука, 1962.--349~с.
\bibitem{HL} Hosmer D. W., Lemeshow S. Applied Logistic Regression. 2nd ed. --- A Wiley-Interscience Publication. 2000.--375 p.
\bibitem{Christensen} Christensen R. Log-Linear Models and Logistic Regression. 2nd ed. --- Springer, 1997.--484 p.
\bibitem{RT} Rao C.R., Toutenburg H. Linear Models: Least Squares and Alternatives. 2nd ed. --- Springer, 1999.--428 p.
\bibitem{KK} Kariya T., Kurata H. Generalized Least Squares. --- Wiley, 2004.--289 p.
\bibitem{LL} Справочник по прикладной статистике. Т. 2. Под ред. Э.Ллойда, У.Ледермана. --- М.: Финансы и статистика, 1990.-- 526 с.
\bibitem{Cochran} Cochran W.G. The omission or addition of an independent variate in multiple linear regression. J. R. Stat.Soc. Suppl.,1938.  \textbf{5}, 171-176.
\bibitem{Quenouille} Quenouille M.H. An application of least squares to family diet surveys. Econometrica,. 1950. \textbf{18}, 27-44.
 \bibitem{MagnusNeu} Магнус Я.\,Р., Нейдеккер Х. Матричное дифференциальное исчисление с
приложениями к статистике и эконометрике. М.: Физматлит. -- 2002. -- 496 с.
 \bibitem{KozlovProchorov} Козлов М.В., Прохоров А.В. Введение в математическую статистику. М: МГУ. 1987. --- 264 с.
\bibitem{PV} Панов В.Г. Вараксин А.Н. О связи между коэффициентами простой и множественной регрессионных моделей. // Сиб. Мат. Ж. 2010. Т. 51, \No 1. С. 196--203.
\bibitem{VMChA}Вараксин А.\,Н., Маслакова Т.\,А., Чуканов В.\,Н., Антонов К.\,Л..
Регрессионная модель зависимости заболеваемости населения от степени загрязнения атмосферного воздуха // Экологические системы и
приборы. 2004. № 4. С. 52--55.
 \bibitem{MVCh} Маслакова Т.\,А., Вараксин А.\,Н., Чуканов В.\,Н. Интерпретация прогностических
регрессионных моделей в области медико-экологического мониторинга // Экологические системы и приборы. 2008. № 2. С. 6--9.
 \bibitem{Marcus-Minc} Маркус М., Минк Х. Обзор по теории матриц и матричных неравенств. Москва.:
Едиториал УРСС, 2004. --- 232 с.
 \bibitem{Gantmakher} Гантмахер Ф.\,Р. Теория матриц. М.: Наука, 1988.--548 с.

\end{thebibliography}
\end{document}